\documentclass[12pt]{article}
\usepackage{times}
\usepackage{mathrsfs}
\usepackage{amssymb,amsmath,amsthm} %
\usepackage{graphicx,psfrag}
\usepackage{amsfonts}
\newcommand{\overbar}[1]{\mkern 1.5mu\overline{\mkern-1.5mu#1\mkern-1.5mu}\mkern 1.5mu}

\def\sp4{\vskip 4mm}
\newcommand \textbox[1]{0.5\textwidth}
\def\@fnsymbol#1{\ensuremath{\ifcase#1\or *\or \dagger\or \ddagger\or
   \mathsection\or \mathparagraph\or \|\or **\or \dagger\dagger
   \or \ddagger\ddagger \else\@ctrerr\fi}}
 
\usepackage{lineno}
\usepackage{mdframed}
\begin{document}
\global\def\refname{{\bf References:}}
%
\baselineskip 26.0pt
%
%
\title{{\LARGE  On the $2$-Vertex Fault Hamiltonicity for Graphs}\\ {\LARGE satisfying Ore's Theorem\footnote{This research was partially supported by the Ministry of Science and Technology of the Republic of China under contract MOST 106-2115-M-033-003.}}}
\date{}
\author{\hspace*{-10pt}
\begin{minipage}[t]{2.7in} \normalsize \baselineskip 18.0pt
\centerline{Hsiu-Chunj Pan$^1$, Hsun Su$^2$, and Shin-Shin Kao$^1$\footnote{Corresponding to: Professor S.-S. Kao, e-mail: shin2kao@gmail.com.}}
\centerline{$^1$\emph{Department of Applied Mathematics,}}
\centerline{\emph{Chung Yuan Christian University,}}
\centerline{\emph{Chung Li City, Taiwan 32023, R.O.C.}}
\centerline{$^2$\emph{Department of Public Finance and Taxation,}}
\centerline{\emph{Takming University of Science and Technology,}}
\centerline{\emph{Taipei, Taiwan 11451, R.O.C.}}
\end{minipage} \kern 0in
\\ \\ \hspace*{-21pt}
\begin{minipage}[b]{5.5in} \normalsize
\baselineskip 22.0pt
\begin{center}
\center{\bf Abstract}\
\end{center}
\indent
\hspace{12 pt}For any undirected and simple graph $G = (V,E)$, where $V$ denotes the vertex set and $E$ the edge set of $G$. $G$ is called hamiltonian if it contains a cycle that visits each vertex of $G$ exactly once. Ore (1960) proved that $G$ is hamiltonian if  $deg_G(u)+deg_G(v)\geq n$ holds for every nonadjacent pair of vertices $u$ and $v$ in $V$, where $n$ is the total number of distinct vertices of $G$.  Kao et al. (2012) proved that any graph $G$ satisfying Ore's theorem remains hamiltonian after the removal of any vertex $x\in V$ unless $G$ belongs to one of the two exceptional families of graphs. In this paper, we proved that in fact, any graph satisfying Ore's theorem remains hamiltonian after the removal of two vertices $x, y \in V$ unless $G$ belongs to one of the nine exceptional families of graphs.
\end{minipage}
\vspace{-10pt}}
\maketitle 
%
%
\clearpage
\section{Introduction}
\label{S1} \vspace{-4pt}
\indent \indent In this paper, we follow the graph definitions and notations of [2], and consider undirected and simple graphs only.  $G=(V,E)$ is a graph if $V$ is a finite set and  $E \subseteq \{(u,v)\mid (u,v)$  is an unordered pair of $V\}$, where $V$ is the vertex set and $E$ is the edge set of $G$. We use $|G|$ or $|V|$ for the number of distinct vertices in $G$,  $K_n$ the complete graph with $n$ vertices, and $\overbar{K_n}$ the graph with $n$ isolated vertices. Two vertices $u$ and $v$ of  $G$ are \emph{adjacent} if $(u,v)\in E$. Given a vertex $u$ of $G$, the \emph{neighborhood} of $u$ is the set  $\{v \mid (u,v)\in E \}\subseteq V$, denoted by   $N_G(u)$. The \emph{degree} of $u$,  $deg_G(u)$, is defined by $deg_G(u)=|N_G(u)|$. The \emph{minimum degree} of $G$, denoted by $\delta(G)$, is min$\{deg_G(u) | u \in V(G)\}$;  the notation $\sigma_2(G)$ is defined by $\sigma_2(G)=min\{deg_G(u)+deg_G(v)|u$ and $v$ are non-adjacent vertices of $G$\}.  Let $S$ be a subgraph of $G$. Define two symbols $N_s(u)= N_G(u)\cap S$, and $deg_s(u)=\mid N_s(u)\mid$.  A \emph{path} in a graph is a single vertex or an ordered list of distinct vertices $v_0, v_1,\cdots ,v_k$ such that $(v_{i-1},v_i)$ is an edge for $1\leq i \leq k$. The first and the last vertices of a path are its \emph{endpoints}.  A cycle is a path of at least three vertices among which the first vertex is the same as the last vertex. A path(cycle) is a hamiltonian path (hamiltonian cycle) if it traverses all vertices of $V$ exactly once.
 A \emph{hamiltonian graph} is a graph with a hamiltonian cycle. A non-hamiltonian graph $G$ is \emph{maximal} if the addition of any edge transforms the graph into a hamiltonian one.  We use $C_n$ for a cycle with $n$ vertices. A graph $G$ is \emph{connected} if it has a path from $u$ to $v$ for each pair of distinct vertices $u, v \in V(G)$.  
 A \emph{vertex cut} of a graph $G$ is a set $S \subseteq V(G)$ such that $G-S$ has more than one component. A graph is \emph{k-connected} if every vertex cut has at least $k$ vertices. The connectivity of $G$, written as $\kappa (G)$, is the minimum size of a vertex cut. In other words, $\kappa (G)$ is the maximum $k$ such that $G$ is $k$-connected.
\\
\indent
The following results about any graph being hamiltonian are well-known.
\vspace{6 pt}\\ 
\noindent \textbf{Theorem 1} (See Ore, 1960 [3]) \emph{A simple graph $G =(V,E) $ with $ | G |=|V| \geq 3$ is hamiltonian if, for each pair of nonadjacent vertices $u$ and $v$ in $V$, $deg_G(u) +deg_G(v) \geq n$.}\\ 
\indent Let $G=(V,E)$ be a graph and let $E' \subseteq E$; we use $G-E'$ to represent the subgraph obtained by removing $E'$ from $G$. Let $F \subseteq V \cup E$. We use $G-F$  for the graph where $V(G-F) =V-F \cap V$ and $E(G-F) = E -\{e| e$ is adjacent to any vertex in $F \cap V\}-E \cap F$.
Suppose that $G-F$ is hamiltonian for any $F\subseteq V\cup E$ and $|F|\leq k$, then $G$ is called a \emph{k-fault-tolerant} \emph{hamiltonian} graph.  
If $F \subseteq V$ and $|F|\leq k$, $G$ is called a $k$-vertex-\emph{fault-tolerant hamiltonian} graph; if $F \subseteq E$ and $|F|\leq k$, $G$ is called a $k$-edge-\emph{fault-tolerant hamiltonian} graph.
It is easy to see that every \emph{k-fault-tolerant (k-vertex-fault-tolerant}, \emph{k-edge fault-tolerant}, respectively) \emph{hamiltonian} graph has at least $k+3$ vertices [2]. 
Moreover, the degree of each vertex in a \emph{k-fault-tolerant} (\emph{k-vertex-fault-tolerant}, \emph{k-edge fault-tolerant}, respectively) \emph{hamiltonian} graph is at least $k+2$ [2].\\ 
\indent We define several operations for graphs.  Let $G_1=(V_1, E_1)$, $G_2=(V_2, E_2)$ be two simple graphs. We say that $G_1$ and $G_2$ are \emph{disjoint} if they have no vertex in common, and \emph{edge-disjoint} if they have no edge in common. The \emph{union} of $G_1$ and $G_2$, denoted by $G_1 \cup G_2$,  is a graph with $V(G_1 \cup G_2)=V(G_1) \cup V(G_2)$ and $E(G_1 \cup G_2)=E(G_1) \cup E(G_2)$; if $G_1 $ and $G_2$ are disjoint, we sometimes denote their union by $G_1+G_2$, and the union of  $k$ copies of $G_1$ by $kG_1$. The \emph{join} of disjoint graphs $G_1 $ and $G_2$, denoted by $G_1 \vee G_2$,  is the graph obtained from $G_1+G_2$ by joining each vertex of $G_1$ to each vertex of $G_2$. Let $H_i$ be any simple graph with $i$ vertices. Define two families of graphs: $\mathscr{G}_1 \equiv K_3 \cup \{H_2 \vee (K_s + K_t) |s + t = n-2, s \geq 1, t \geq 1\}$, and $\mathscr{G}_2\equiv \{H_s \vee sK_1  | 2s = n\}$. See Fig. 1.\\
 \noindent \textbf{Theorem 2} (See Su, Shih, and Kao, 2012 [6])  \emph{Let $G = (V, E)$ be a graph with $|G|=|V|=n\geq 3$. Suppose that $deg_G(u) + deg_G(v) \geq n$ holds for any nonadjacent pair $\{u, v\} \subset V$, then either $G$ is 1-vertex-fault hamiltonian or $G$ belongs to one of the two families $\mathscr{G}_1$ and $\mathscr{G}_2$. 
In addition, $G$ is either 1-edge-fault hamiltonian or  $G \in  \mathscr{G}_1$ with  $s \in \{1, 2\}$.}
\vspace{5 pt}
\vfill
\begin{figure}[h!]
\graphicspath{{graphs/}}
\begin{center}
\centering
\includegraphics[scale=1.45]{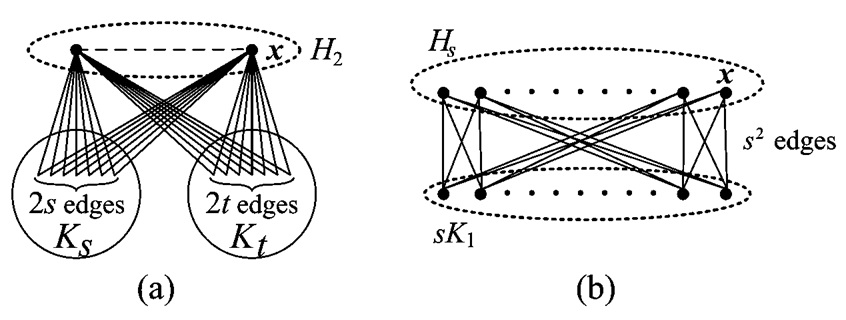}
\vspace{-15pt}
\end{center}
\center{\text{Figure 1:  An illustration of graphs of (a) $\{H_2 \vee (K_s + K_t) |s + t = n - 2$, $s \geq 1$, $t \geq 1\}$ in $\mathscr{G}_1$; (b) $\mathscr{G}_2$}}.
\vspace{5pt}
\end{figure}
\indent In this paper, we want to show that any graph satisfying the degree-sum condition in Theorem 1 is not only a hamiltonian graph but also a \emph{2-vertex-fault hamiltonian} graph, unless it belongs to certain exceptional families.  It is clearly that $G$ is not 2-vertex-fault tolerant (resp. 2-edge-fault tolerant) when the vertex-connectivity (resp. edge-connectivity) of a graph $G$ is equal to or less than 3. Thus, in the Theorem 6 of this paper we only consider graphs whose vertex connectivity is greater than or equal to 4.
\section{Main Results}
\indent \indent Following the discussion in the previous section, we obtain the graph $G_1\colon G_2$ from $G_1+G_2$ by connecting some vertices of $G_1$ to some vertices of $G_2$, possibly with constraints on how edges are added. For three simple graphs $G_1$, $G_2$ and $G_3$, the notation $G_1\colon G_2\colon G_3$ is defined by $G_1\colon G_2\colon G_3$=$(G_1\colon G_2)\colon G_3$. So $G_1\colon G_2\colon G_3$ is the graph obtained from $G_1+G_2+G_3$ by connecting some vertices of $G_s$ to some vertices of $G_t$, where $s,t \in \{1,2,3\}$ and $s\ne t$; there may be restrictions on how edges are added. Thus, $G_1+G_2 \subseteq G_1\colon G_2 \subseteq G_1\vee G_2$, and $(G_1+G_2+G_3)\subseteq G_1\colon G_2\colon G_3 \subseteq (G_1\vee G_2) \vee G_3$. An example for $G_1\colon G_2\colon G_3$ is given by $H_i\colon x\colon y$, where $H_i$ is a simple graph with $i$ vertices, and $x$ and $y$ are two vertices not belonging to $V(H_i)$.  Then  $H_i+x+y \subseteq H_i\colon x\colon y \subseteq (H_i\vee x)\vee y$. The \emph{join-and-delete-one-edge} of $G_1\vee G_2$, written as $G_1 \vee^-  G_2$, is obtained by deleting one edge from the graph $G_1 \vee G_2$. 
The \textit{join-and-delete-two-edges} of $G_1\vee G_2$, written as $G_1 \vee^= G_2$, is obtained by deleting two edges from the graph $G_1\vee G_2$. Both $G_1 \vee^- G_2$ and $G_1 \vee^= G_2$ are edge-deleted subgraphs of $G_1\vee G_2$[1]. \\ 
%
\indent We define seven graph families in the following.\\
\noindent \textbf{Definition 3} Let $H_k$ be any simple graph with $k$  vertices.  Define $\eta_i$ for $1\le i\le 7$ as below. \\
(1) $\eta_1 =(H_{(n+1)/2}\vee \overbar{K_{(n-1)/2}}) - (v_1, u_1)$,  $v_1\in V(H_{(n+1)/2})$, $u_1\in V(\overbar{K_{(n-1)/2}})$,
$n$ is odd, and $n\geq 9$, with $ deg_{H_{(n+1)/2}}(v_1)\ge 2$.  and $\sigma_2(H_{(n+1)/2}) \geq1 $.  See Figure 2.\\
(2) $\eta_2= (H_{(n+1)/2} \vee \overbar{K_{(n-1)/2}})$,  $n$ is odd, $n \geq 7$. When $n = 7$, $\sigma_2(H_4) > 1$, and $\kappa(H_4) \ge 1$; when $n\geq 9$,  $\sigma_2(H_{(n+1)/2})\ge 1 $. See Figure 3 and Figure 4.\\ 
(3) $\eta_3= H_{n/2} \vee (\overbar{K_{(n-4)/2}} \cup K_2)$, $n$ is even, and $n \ge 8$. See Figure 5.\\
(4) $\eta_4= (H_{(n-4)/2} \colon K_2)\vee (\overbar{K_{(n-4)/2}}\cup K_2)$, $n$ is even, and $n \ge 8$.  See Figure 6.\\
   \indent \hspace{-5pt}$\eta_4^- = (H_{(n-4)/2} \colon K_2)\vee ^- (\overbar{K_{(n-4)/2}}\cup K_2)$,  $n$ is even, and $n \ge 8$.  See Figure 7.\\
   \indent  \hspace{-5pt}$\eta_4^= = (H_{(n-4)/2} \colon K_2)\vee ^= (\overbar{K_{(n-4)/2}}\cup K_2)$, $n$ is even, and $n \ge 8$.  See Figure 8.\\
\indent In the graph $\eta_4=(H_{(n-4)/2}\colon K_2)\vee  (\overbar{K_{(n-4)/2}}\cup K_2)$, $V(H_{(n-4)/2}\colon K_2) = \{v_1, v_2, \cdots , v_{(n-4)/2}\}\cup \{x, y\}$, $V(\overbar{K_{(n-4)/2}}\cup K_2)= \{u_3, u_4,\cdots , u_{n/2}\}\cup \{u_1, u_2\}$. Let $RE=\{(x, u_1), (x, u_2), (y, u_1), (y, u_2)\}$. The graph  $\eta_4^-$ is obtained by deleting one edge $e\in RE$ from  $\eta_4$. The graph  $\eta_4^=$ is obtained by deleting two vertex-disjoint edges $e_1, e_2\in RE$ from $\eta_4$. \\
(5) $\eta_5= H_{n/2} \vee \overbar{K_{n/2}}$, $n$ is even, $n \ge 8$.  See Figure 9.\\
(6) $\eta_6= H_4 \vee 3K_2$.  See Figure 10.\\
(7) $\eta_7=  H_4 \vee (2K_2 \cup K_1)$.  See Figure 11.\\ 
\vfill
\begin{figure}[h! bt]
\graphicspath{{graphs/}}
\begin{center}
\centering
\includegraphics[scale=0.87]{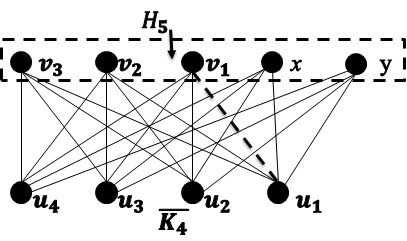}
\vspace{-10pt}
\end{center}
\center{\parbox{\textwidth}{Figure 2:  $Let$ $n = 9$, $\eta_1 = (H_5\vee \overbar{K_4})-(v_1, u_1)$, $with$ $\sigma_2(H_5)\ge 1$, $deg_{H_5 }(v_1)\ge 2$,\\
$where$ $V(H_5)=\{v_1,v_2,v_3,x,y\}$, $V(\overbar{K_4})=\{u_1,u_2,u_3,u_4\}$. $Obviously$, $deg_{\eta_1}(x)\ge 4$,  $deg_{\eta_1}(y)\geq 4$, $deg_{\eta_1}⁡(v_i )\geq 4$,  for $i=1, 2, 3$. $deg_{\eta_1}⁡(u_i )=5$ for $i=2,3,4$. $deg_{\eta_1}⁡(u_1 )=4$.}}\\
\vspace{5pt}
\vspace{5pt}
\graphicspath{{graphs/}}
\centering
\includegraphics[scale=0.87]{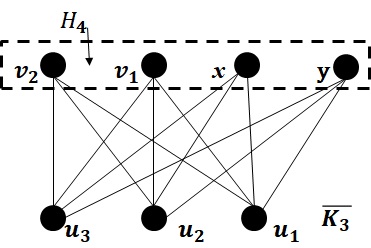}
\vspace{-2pt}
\center{\parbox{\textwidth}{Figure 3: $Let$ $n=7$, $\eta_2 = (H_4\vee \overbar{K_3})$, $with$ $\sigma_2(H_4) > 1$, and $\kappa(H_4)\ge 1$. $V(H_4)=\{v_1,v_2,x,y\}$, $V(\overbar{K_3})=\{u_1,u_2,u_3\}$. Then $deg_{\eta_2}(v_i) > 3$ for $i=1, 2$; $deg_{\eta_2}(x)> 3$, $deg_{\eta_2}(y)> 3$, and $deg_{\eta_1}⁡(u_i )=4$ for $i=1, 2, 3$.}} 
\vspace{5pt}
\end{figure}
\clearpage
\vspace{5 pt}
\begin{figure}[!h]
\graphicspath{{graphs/}}
\begin{center}
\centering
\includegraphics[scale=0.9]{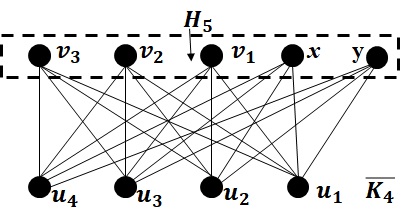}
\vspace{-2pt}
\end{center}
\center{\parbox{\textwidth}{Figure 4: $Let$ $n = 9$, $\eta_2 = (H_5\vee \overbar{K_4}) $, with $\sigma_2(H_5)\ge 1$, $V(H_5)=\{v_1,v_2,v_3,x,y\}$, $V(\overbar{K_4})=\{u_1,u_2,u_3,u_4\}$. Obviously,  $deg_{\eta_2}⁡(v_i )\geq 4$ for $i=1, 2, 3$. $deg_{\eta_2}(x)\ge 4$, $deg_{\eta_2}(y)\geq 4$, and $deg_{\eta_1}⁡(u_i )=5$ for $i=1, 2, 3, 4$.}}
\vspace{15pt}
%
\graphicspath{{graphs/}}
\begin{center}
\centering
\includegraphics[scale=0.7]{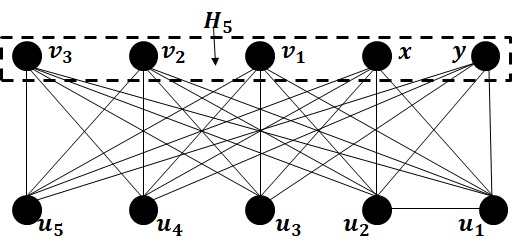}
\vspace{-2pt}
\end{center}
\center{\parbox{\textwidth}{Figure 5: $Let$ $n=10$, $\eta_3= H_5\vee (\overbar{K_3}\cup K_2)$, where $V(H_5)=\{v_1,v_2,v_3, x, y \}$, $V (\overbar{K_3}) ̅=\{u_3, u_4, u_5\}$, $V(K_2)=\{u_1, u_2\}.$ $deg_{\eta_3}(v_i)\geq 5$ for $i=1, 2, 3$; $deg_{\eta_3}(x)\geq 5$, $deg_{\eta_3}(y)\geq 5$ and $deg_{\eta_3}(u_i)=5$ for $i\neq 1, 2$; $deg_{\eta_3}(u_i)=6$ for $i=1, 2$.}}
\vspace{25pt}
%
\graphicspath{{graphs/}}
\begin{center}
\centering
\includegraphics[scale=0.7]{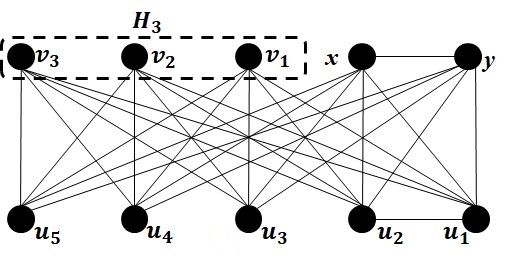}
\vspace{-0pt}
\end{center}
\center{\parbox{\textwidth}{Figure 6: $Let$ $n=10$, $\eta_4= (H_3\colon K_2)\vee (\overbar{K_3}\cup K_2)$. Then $V(H_3)=\{v_1, v_2, v_3\}$, $V(\overbar{K_3})=\{u_3, u_4, u_5\}$, the complete graph in $(H_3\colon K_2)$ is with $V(K_2)=\{x, y\}$, and the complete graph in $(\overbar{K_3}\cup K_2)$ is with $V(K_2)=\{u_1, u_2\}$; $deg_{\eta_4}(v_i)\ge 5$ for $i=1,2,3$. $deg_{\eta_4}⁡(x)\ge 6$, $deg_{\eta_4}⁡(y)\ge 6$, and $deg_{\eta_4}⁡(u_i)=5$ for $i\ne 1,2,$; $deg_{\eta_4}⁡(u_i)=6$ for $i=1, 2 $.}} 
\vspace{15pt}
\end{figure}
\clearpage
\begin{figure}[!ht]
\graphicspath{{graphs/}}
\begin{center}
\centering
\includegraphics[scale=0.7]{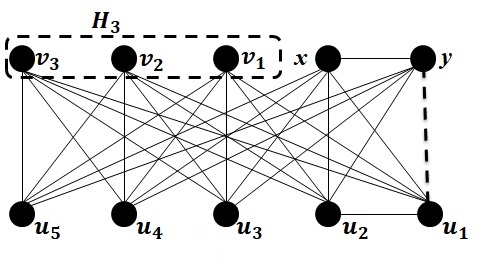}
\vspace{3pt}
\end{center}
\center{\parbox{\textwidth}{Figure 7: $Let$ $n=10$, $\eta_4^- = (H_3\colon K_2)\vee^-  (\overbar{K_3}\cup K_2)$.  Then $V(H_3)=\{v_1,v_2,v_3\}$, $V(\overbar{K_3})=\{u_3, u_4, u_5\}$, the complete graph in $(H_3\colon K_2)$ is with  $V(K_2)=\{x, y\}$, and the complete graph in $(\overbar{K_3}\cup K_2)$ is with $V(K_2)=\{u_1, u_2\}$; $deg_{\eta_4^- }(v_i)\ge 5$ for $i=1,2,3$; $deg_{\eta_4^- }(x)\ge 6$, $deg_{\eta_4^- }(y)\ge 5$; $deg_{\eta_4^- }⁡(u_i )=5$ for $i\ne 2$; $deg_{\eta_4^- }(u_2)=6$.}}
\vspace{5pt}
%
\graphicspath{{graphs/}}
\begin{center}
\centering
\includegraphics[scale=0.7]{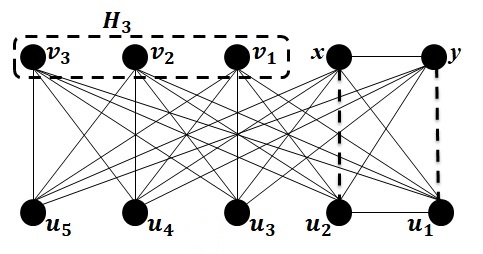}
\vspace{5pt}
\center{\parbox{\textwidth}{Figure 8: $Let$ $n=10$, $\eta_4^== (H_3\colon K_2)\vee^= (\overbar{K_3}\cup K_2)$.  Then  $V(H_3)=\{v_1,v_2,v_3\}$, $V(\overbar{K_3})=\{u_3, u_4, u_5\}$, the complete graph in $(H_3\colon K_2)$ is with  $V(K_2)=\{x, y\}$, and the complete graph in $(\overbar{K_3}\cup K_2)$ is with $V(K_2)=\{u_1, u_2\}$; $deg_{\eta_4^=} (v_i)\ge 5$ for  $i=1, 2, 3$; $deg_{\eta_4^=}(x)\ge 5$, $deg_{\eta_4^=}(y)\ge 5$, $deg_{\eta_4^=}(u_i)=5$ for  $i=1, 2, 3, 4, 5$.}}
\vspace{5pt}
\end{center}
%
\graphicspath{{graphs/}}
\begin{center}
\centering
\includegraphics[scale=0.7]{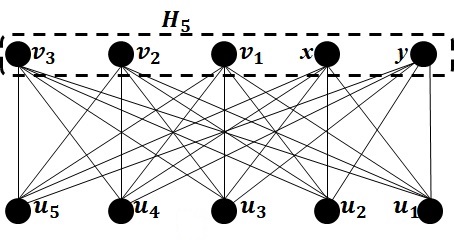}
\vspace{-5pt}
\end{center}
\center{\parbox{\textwidth}{Figure 9: $Let$ $n=10$, $\eta_5= H_5\vee \overbar{K_5}$. Then $V(H_5)=\{v_1,v_2,v_3, x, y\}$, $V(\overbar{K_5})=\{u_1, u_2, u_3, u_4, u_5\}$; $deg_{\eta_5}(v_i)\geq 5$ for $i=1,2,3$; $deg_{\eta_5}(x)\geq 5$, $deg_{\eta_5}(y)\geq 5$, $deg_{\eta_5}(u_i)=5$ for $i=1,2,3,4,5$.}}
\vspace{5pt}
\end{figure}
\clearpage
\begin{figure}[!ht]
\graphicspath{{graphs/}}
\begin{center}
\centering
\includegraphics[scale=0.7]{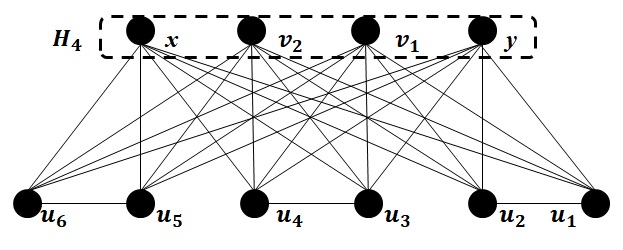}
\vspace{-3pt}
\end{center}
\center{\parbox{\textwidth}{Figure 10:  $\eta_6= H_4\vee 3K_2$. Then $V(H_4)=\{v_1,v_2, x, y\}$, The three complete graphs in $3K_2$ are with  $V(K_2)=\{u_i, u_{i+1}\}$, for $i=1,3,5$. $deg_{\eta_6}(v_i)\ge 6$ for $i=1, 2$. $deg_{\eta_6}(x)\ge 6$, $deg_{\eta_6}(y)\ge 6$, and $deg_{\eta_6}(u_i)=5$ for $i=1, 2, 3, 4, 5, 6$.}}
\vspace{10pt}
%
\graphicspath{{graphs/}}
\begin{center}
\centering
\includegraphics[scale=0.7]{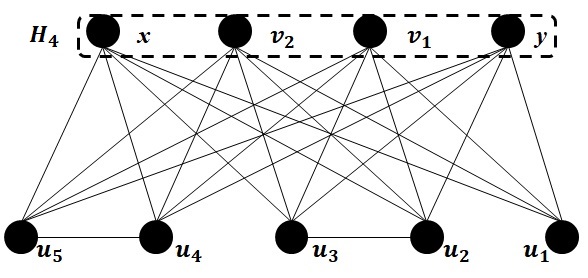}
\vspace{-2pt}
\end{center}
\center{\parbox{\textwidth}{Figure 11:  $\eta_7= H_4\vee (2K_2\cup K_1)$. Then $V(H_4)=\{v_1,v_2,x, y\}$. The two complete graphs in $2K_2$ are with  $V(K_2)=\{ u_i, u_{i+1}\}$, for $i=2, 4$; and $V(K_1)=\{u_1\}$. $deg_{\eta_7}(v_i)\ge 5$ for $i=1, 2$. $deg_{\eta_7}(x)\ge 5$, $deg_{\eta_7} (y)\ge 5$, $deg_{\eta_7}⁡(u_i )=5$ for $i\ne 1$ and  $deg_{\eta_7}⁡(u_1 )=4.$ }}\\
\end{figure}
\vspace{0pt}
\vfill
\indent In the sequel, we let $G' =(V' , E' )$ be any graph with $|V'|=n' \geq 3 $ such that for any nonadjacent vertices $u$ and $v \in V'$, $deg_{G' }(u)+deg_{G' }(v) \geq n' $. Let $G =(V, E)$ with $V= V' \cup \{x\}$.  We consider the following two cases for $G$:\\  
(1) If $ E = E' \cup  \{(x,y)| y = u, v$,  where $deg_{G' }(u) + deg_{G' }(v) \geq |G'|$ and  $u, v$ are nonadjacent in $G' \}$, then for any nonadjacent vertices $u$ and $v$, $deg_G(u) + deg_G(v)\geq |G|+1$.\\
(2) If $E = E'\cup (\{(x,y)| y = u,v$, where $deg_{G'}(u)+deg_{G'}(v) \geq |G'| $ and  $u, v$ are nonadjacent in $G'$ $\}-e)$. Then for any nonadjacent vertices $u$ and $v$, $deg_G(u)+deg_G(v) \geq |G|$.\\
\indent In 1985, Ainouche and Christofides gave the following result.\\
\noindent \textbf{Theorem 4}  (See Ainouche and Christofides, 1985 [4]) \emph{Let $G=(V, E)$ be a 2-connected maximal non-hamiltonian graph. If $deg_G(a) +deg_G(b) \geq |G|-2$ for any two non-adjacent vertices $a$, $b$, then $G$ is isomorphic to one of the following five graphs: $G_1= K_{(n-1)/2}\vee \overbar{K_{(n+1)/2}}$, $n$ is odd, $n \ge 3$; $G_2=K_{(n-2)/2}\vee \overbar{K_{(n+2)/2}}$, $n$ is even, $n \ge 4$; $G_3= K_{(n-2)/2} \vee (\overbar{K_{(n-2)/2}}\cup K_2)$, $n$ is even, $n \ge 4$; $G_4= K_2 \vee (2 K_2 \cup K_1)$; $G_5= K_2 \vee 3 K_2$.} \\
\indent In 2013, Zhao presented the following result:  \\
\noindent \textbf{Theorem 5} (See Zhao, 2013 [5]) \emph{If $G''$ is a connected graph of order $n \geq 3$ such that $deg_{G''}(x)+deg_{G''}(y) \geq n-2 $ for each pair of nonadjacent vertices $x, y$ in $G''$, then either $G''$ is hamiltonian or $G''$ is isomorphic to one of the following nine graphs: }\\
\vspace{10pt}
\noindent \hspace{-5pt}\emph{(1)$K_{1, 3}$, (2)$H_2 \vee 3K_2 $, (3)$H_2 \vee (2K_2 \cup K_1) $, (4)$K_h \colon \omega \colon K'_t$, (5)$(H_{(n-1)/2} \vee \overbar{K_{(n+1)/2}})-e$ , (6)$K_1 \colon C'_6 $, (7)$H_{(n-1)/2} \vee \overbar{K_{(n+1)/2}}$ ,  (8)$H_{(n-2)/2} \vee \overbar{K_{(n+2)/2}}$ , (9)$H_{(n-2)/2} \vee (\overbar{K_{(n-2)/2}} \cup K_2)$.}\\
\indent The graphs in (4) and (6) given in Theorem 5 above need more explanations. 
For $G'' \in K_h  \colon w \colon K'_t$ in (4), $K'_t$ is a graph by removing some (none, one, or more) vertex-disjoint edges of $K_t$, with $h \le t$; in the operation of ``:'', edges are added from $w$ to $K_h$ and $K'_t$ as long as $\sigma_2(G'')\ge |G''|-2$ holds.
It is easy to see that $G''-\{\omega\}$ is a disconnected graph. For the case of $G''=K_1\colon C'_6$ in (6), see Fig. 12. It can be seen that $G''=\overbar{K_3}\vee \overbar{K_4}-(u, x_5)$. Thus $G''$ belongs to $H_3\vee \overbar{K_4}-e$, which is in (5) of Theorem 5.\\ 
\vfill
\begin{figure} [!h]
\graphicspath{{graphs/}}
\begin{center}
\centering
\includegraphics[scale=1.0]{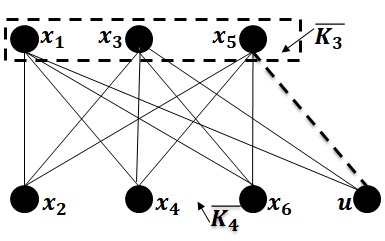}
\vspace{-15pt}
\end{center}
\center{Figure 12:  $G''=(K_1\colon C'_6)=\overbar{K_3}\vee \overbar{K_4}-(u, x_5)$.}
\end{figure}
\vspace{12pt}
\clearpage\clearpage
\begin{figure} [!ht]
\graphicspath{{graphs/}}
\begin{center}
\centering
\includegraphics[scale=1.0]{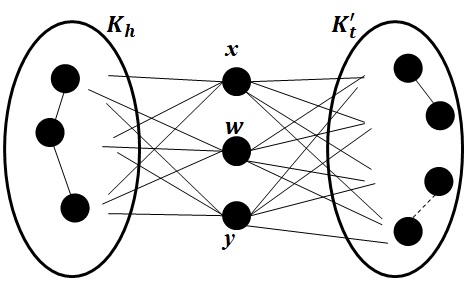}
\vspace{-15pt}
\end{center}
\center{\text{Figure 13:  $G \in (K_h \colon H_3\colon K'_t)$}}
\vspace{13pt}
%
\graphicspath{{graphs/}}
\begin{center}
\centering
\includegraphics[scale=1.0]{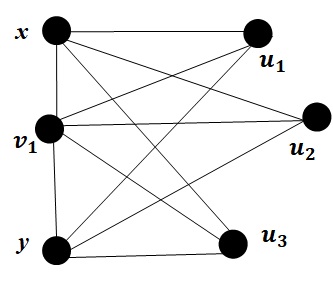}
\vspace{-15pt}
\end{center}
\center{\text{Figure 14:  $G=H_3\vee \overbar{K_3}$}}
\vspace{2pt}
\end{figure}
\vspace{1pt}
\noindent 
\textbf{Theorem 6} \emph{Let $G =(V,E)$ be a graph with $|G|=|V|=n \geq 5$ and $\kappa (G)\geq 4$. Suppose that  $ \sigma_2(G)\geq n$ holds for any nonadjacent pair $\{u, v\}\subset V$, then either $G$ is 2-vertex-fault hamiltonian or
$G \in \{\eta_1, \eta_2, \eta_3, \eta_4, \eta_4^-, \eta_4^=, \eta_5, \eta_6, \eta_7\}$}.\\
\begin{proof}
Given a graph  $G =(V,E)$ satisfying the conditions of the theorem, we can obtain $G''$ from $G$ by deleting two vertices from $V(G)$, then $| G''|= n''=n-2$, $\kappa (G'')\geq 2$,  and   $\sigma_2(G'')\geq n-4\geq n''-2$.\\
\textbf{Case 1. }  $G''$ is hamiltonian. Thus $G$ is 2-vertex-fault hamiltonian.\\
\textbf{Case 2. }  $G''$ is not hamiltonian. Thus $G''$ belongs to one of the eight graphs in Theorem 5. (Note that $K_1\colon C'_6$ belongs to (5) ).  we obtain $G$ by adding two vertices $x,y$ to $G''$ such that $\sigma_2(G)\geq n$ and $\kappa(G)\ge 4$. From Case 2.1 to Case 2.9, we will show that only some of $G''$can be transformed to specific $G$ by adding vertices or deleting edges.  It can be verified that $\kappa(G)\ge 4$ for each $\eta_i$, except for $\eta_2$ when $n=7$. We need some extra condition so as to get $\kappa(G)\ge 4$,  the reason for which is explained in Case 2.7.\\
\noindent
\textbf{Case 2.1. }Let $G''=K_{1,3}=K_1\vee \overbar{K_3}$. We add two vertices $x$ and $y$ to $G''$ such that $G=H_3\vee \overbar{K_3}$. Since the degree of each vertex $u_i$ is 3 for $u_i\in V(\overbar{K_3})$, the graph $G$ is not 2-vertex-fault tolerant graph. See Fig. 14. Obviously, $H_3\vee \overbar{K_3}$ is a special case for $\eta_5$.\\
\textbf{Case 2.2. }Let $G''=H_2  \vee  3K_2$ , then we have  $\delta(G'') =3$, and $\sigma_2(G'') = 6$. We add two vertex $x$ and  $y$ to $G''$ such that $G =H_4  \vee 3K_2$. It can be seen that $\delta(G) = 5$, $\sigma_2 (G)= 10$. We use $\eta_6$ for this kind of graph,  $\eta_6 = H_4  \vee 3K_2$. See Fig. 10.\\
\textbf{Case 2.3. }Let $G''=H_2 \vee (2K_2 \cup K_1)$. We have $\delta(G'') =2$, and $\sigma_2 (G'')= 5$.
We add two vertices $x$ and $y$ to $G''$ such that $G = H_4 \vee (2K_2\cup K_1)$. It can be seen that $\delta(G) = 4$, $\sigma_2 (G)= 9$. We use $\eta_7$  for this kind of graph, where $\eta_7 = H_4 \vee (2K_2 \cup K_1)$. See Fig. 11.\\
\textbf{Case 2.4. }Let $G'' \in K_h \colon w \colon K'_t$.  When adding two vertices to $G''$, we obtain 
 $G\in K_h \colon H_3  \colon K'_t$, where $\kappa(G) = 3$. See Fig. 13. Consequently, $G$ is not 2-vertex-fault tolerant.\\
\textbf{Case 2.5. }Let $G''=(H_{(n''-1)/2} \vee \overbar{K_{(n''+1)/2}}) - e$, $V(H_{(n''-1)/2}) = \{v_1,v_2,\cdots, v_{(n''-1)/2}\}$, and $V(\overbar{K_{(n''+1)/2}}) = \{u_1,u_2,\cdots , u_{(n''+1)/2}\}$,$e=(v_1, u_1)$. Obviously $\delta (G'')=(n''-1)/2-1$ and $\sigma_2(G'')=n''-2$.
We add two vertices $x$ and $y$ to $G''$ such that 
$G= (H_{(n+1)/2} \vee \overbar{K_{(n-1)/2}})-(v_1, u_1)$, where $V(H_{(n+1)/2}) = \{v_1,v_2, \cdots , v_{(n-3)/2}, x , y \}$,  $V(\overbar{K_{(n-1)/2}}) = \{u_1, u_2, \cdots , u_{(n-1)/2}\}$, $deg_{H_{(n+1)/2}}(v_1)\ge 2$, and $\sigma_2(H_{(n+1)/2})\ge 1$. Note that if $|G|=n \le 7$, there will be a vertex with degree less than four. Thus $G$ is not 2-vertex fault hamiltonian. So we examine $G$ for  $n \ge 9$. 
%
If $n\geq 9$, the following explains why graph $G$  satifies $\sigma_2(G) \geq n$. For graph $H_{(n+1)/2}$ itself, before the operation of ``join'', we have ``$\sigma_2(H_{(n+1)/2})\geq 1$''. 
Once the operation ``join'' is completed, and the edge $(v_1, u_1)$ has been deleted, in graph $G$, except vertex $v_1$, the degrees of those vertices of $V(H_{(n+1)/2})$ are all increased by $(n-1)/2$, and the degree of $v_1$  is increased by $[(n-1)/2 -1]$. Clearly $deg_G(v_1) \ge 2+[(n-1)/2 -1]= (n+1)/2$. At the same time, except vertex $u_1$, the degrees of those vertices of $\overbar{V(K_{(n-1)/2}})$ are all increased by $(n+1)/2$, and $deg_G(u_1)=(n-1)/2$. For any two distinct vertices $a$ and $b$, we have the following four cases.\\
\textbf{Subcase 2.5.1.} Suppose both $a$ and $b$ belong to $H_{(n+1)/2}$, $a\neq v_1$, and $b \neq v_1$ then $deg_G (a)+deg_G⁡(b) = deg_{H_{(n+1)/2}} (a)+deg_{H_{(n+1)/2} }⁡(b) + (n-1)/2 + (n-1)/2\geq 1+ (n-1)=n$. \\
\textbf{Subcase 2.5.2.} Suppose both $a$ and $b$ belong to $H_{(n+1)/2}$, and one of $a$ and $b$ is equal to $v_1$. For convenience, let $a=v_1$, then $deg_G (v_1)+ deg_G⁡(b)\geq 1+(n+1)/2 + (n-1)/2= n+1$.\\
\textbf{Subcase 2.5.3.} Suppose both $a$ and $b$ belong to $\overbar{K_{(n-1)/2}}$ , then $deg_G (a)+deg_G⁡(b) \ge (n+1)/2 + (n-1)/2 =n$. \\
\textbf{Subcase 2.5.4.}  $deg_G (v_1)+deg_G⁡(u_1) \ge (n+1)/2 + (n-1)/2 =n$.\\
From \textbf{Subcase 2.5.1} to \textbf{Subcase 2.5.4}, we can see that $\sigma_2(G)\geq n$. 
Note $\eta_1 = (H_{(n+1)/2} \vee \overbar{K_{(n-1)/2}})- (v_1, u_1)$ ,  $v_1\in H_{(n+1)/2} , u_1\in \overbar{K_{(n-1)/2}}$ , $n$ is odd, $n \geq 9$ , with $deg_{H_{(n+1)/2} }(v_1)\ge 2$ and  $\sigma_2(H_{(n+1)/2}) \geq 1$ . See Fig. 2.\\
%
\textbf{Case 2.6. }This case belongs to Case 2.5.\\
\noindent
\textbf{Case 2.7.} As discussed in Case 2.5, let $G''= (H_{(n''-1)/2} \vee \overbar{K_{(n''+1)/2}}), V(H_{(n''-1)/2}) =\{v_1,v_2,\cdots, v_{(n''-1)/2}\}$,  and  $V(\overbar{K_{(n''+1)/2}}) = \{u_1,u_2,\cdots, u_{(n''+1)/2}\}$. It can be seen that $\delta(G'') = (n''-1)/2$, and $\sigma_2(G'') = n''-1.$ $G$ is defined by adding two vertices $x$ and $y$ to $G''$, and $G= (H_{(n+1)/2}\vee \overbar{K_{(n-1)/2}})$. If $n=7$, $G=(H_4 \vee \overbar{K_3})$, where $V(H_4)= \{v_1,v_2, x , y \}$ and $V(\overbar{K_3})=\{u_1,u_2, u_3 \}$,  $\sigma_2(H_4)>1$, the condition $\kappa(H_4)\ge 1$ is required to make sure $\kappa(G)\ge 4$. See Fig. 3. If $n\geq  9$, $G=(H_{(n+1)/2}\vee (\overbar{K_{(n-1)/2}})$, where $V(H_{(n+1)/2})=\{v_1,v_2,\cdots , v_{(n-3)/2}, x, y\}$; and $V(\overbar{K_{(n-1)/2}})=\{u_1,u_2,\cdots , u_{(n-1)/2}\}$, $\sigma_2(H_{(n+1)/2})\geq 1$. See Fig. 4. We use $\eta_2$ to stand for this kind of graph, $\eta_2=(H_{(n+1)/2}\vee (\overbar{K_{(n-1)/2}})$, $n$ is odd, $n \geq 9$, with $\sigma_2(H_{(n+1)/2}) \geq 1$; when $n=7$, $\eta_2=(H_4\vee \overbar{K_3})$ with $\sigma_2(H_4)>1$ and $\kappa(H_4)\ge 1$. \\
%
\textbf{Case 2.8. }Let $G''= H_{(n''-2)/2}\vee \overbar{K_{(n''+2)/2}}, V(H_{(n''-2)/2}) = \{v_1,v_2,\cdots , v_{(n''-2)/2}\}$, and \\
$V(\overbar{K_{(n''+2)/2}}) = \{u_1,u_2,\cdots , u_{(n''+2)/2}\}$. It can be seen that $\delta(G'') = (n''-2)/2$, and $\sigma_2 (G'') = n''-2$. By adding one vertex $x$ to $G''$ such that $G'=((H_{(n' -3)/2} \colon x) \vee \overbar{K_{(n'+1)/2}})$, we have $\delta(G') = (n'-1)/2$, and $\sigma_2 (G') =n'-1$.  
Then adding one more vertex $y$ to $G'$ such that $G=(H_{(n-4)/2}\colon x\colon y) \vee \overbar{K_{n/2}} =H_{n/2} \vee \overbar{K_{n/2}}$, we have $\delta(G)=(n/2)$, and $\sigma_2(G)=n$.  Note that $\eta_5 = H_{n/2}\vee \overbar{K_{n/2}}$. It can be seen that $n \ge 8$ is required for ensuring $\delta (G)\ge 4$. See Fig. 9. It can be observed that $\eta_5$ is isomorphic to $\mathscr{G}_2$.\\
\textbf{Case 2.9. }Let $G''= H_{(n''-2)/2} \vee (\overbar{K_{(n''-2)/2}} \cup K_2)$, $V(H_{(n''-2)/2}) =\{v_1,v_2,\cdots , v_{(n''-2)/2}\}$, \\
$V(K_2)=\{ u_1,u_2\} $ and $V(\overbar{K_{(n''-2)/2}}) = \{u_3,u_4,\cdots , u_{(n''+2)/2}\}$. It can be seen that $\delta(G'') = (n''-2)/2$, and $\sigma_2(G'') = n''-2$. \\
\textbf{Subcase 2.9.1. }Adding one vertex $x$ to $G''$ such that
$G'=((H_{(n'-3)/2}\colon x) \vee \overbar{K_{(n'-3)/2}}\cup K_2)$.
we have  $\delta(G') = (n'-1)/2$, and $\sigma_2(G') = n'-1$. 
Then adding one more vertex $y$ to $G'$ such that 
$G = ((H_{(n-4)/2}\colon x\colon y) \vee(\overbar{K_{(n-4)/2}}\cup K_2)) =H_{n/2} \vee (\overbar{K_{(n-4)/2}}\cup K_2)$,  we have $\delta(G) = (n/2)$, and $\sigma_2 (G) = n$.  It is noted that $\eta_3=H_{n/2}\vee (\overbar{K_{(n-4)/2}}\cup K_2)$.  This Subcase also shows that $n \ge 8$ is required to make sure $\delta (G)\ge 4$.  See Fig. 5.\\
\textbf{Subcase 2.9.2. }Adding one more vertex $y$ and an edge $(x, y)$ to $G'$, where $G'=(H_{(n' -3)/2}\colon  x) \vee (\overbar{K_{(n'-3)/2}} \cup K_2)$ as in Subcase 2.9.1, we obtain the graph 
$G= ((H_{(n-4)/2}\colon x \colon y) \vee (\overbar{K_{(n-4)/2}}\cup K_2)) \cup (x, y)= (H_{(n-4)/2}\colon K_2) \vee (\overbar{K_{(n-4)/2}}\cup K_2)$, which is defined to be $\eta_4$. Let $V(H_{(n-4)/2}\colon K_2)=\{v_1, v_2, \cdots , v_{(n-4)/2}\}\cup \{x, y\}$, $V(\overbar{K_{(n-4)/2}}\cup K_2)= \{u_3, u_4,\cdots , u_{n/2}\}\cup \{u_1, u_2\}$.
We have $\sigma_2(G) = n$. See Fig. 6.\\
\indent Let $RE=\{(x, u_1), (x, u_2) , (y, u_1) , (y, u_2) \}$. For the graph $\eta_4$, since $deg_{\eta_4}⁡(u_1)=deg_{\eta_4}⁡(u_2)=(n/2)+1$, 
$deg_{\eta_4}⁡(x)\ge (n/2)+1$, and $deg_{\eta_4}⁡(y)\ge (n/2)+1$, we can reduce the degree of $u_1$ or $u_2$ by removing one edge(or two edges) that belongs (belong) to $RE$ from $\eta_4$. The following discussions show that to make sure $\delta (G)\ge 4$, $n \ge 8$ is required. \\
\textbf{Subcase 2.9.2.1. }Removing one edge $e \in RE$ from $\eta_4$, we have $\eta_4^-= (H_{(n-4)/2}\colon K_2)\vee^- (\overbar{K_{(n-4)/2}}\cup K_2)=(H_{(n-4)/2}\colon K_2)\vee  (\overbar{K_{(n-4)/2}}\cup K_2)- e$, which shows that $\delta (\eta_4^-)=(n/2)$, and $\sigma_2(\eta_4^-) =n $. See Fig. 7. \\
\textbf{Subcase 2.9.2.2. }Removing two vertex-disjoint edges $e_1, e_2 \in RE$ from $\eta_4$, we have \\
$\eta_4^== (H_{(n-4)/2}\colon K_2)\vee^= (\overbar{K_{(n-4)/2}}\cup K_2)=(H_{(n-4)/2}\colon K_2)\vee  (\overbar{K_{(n-4)/2}}\cup K_2)- e_1-e_2$.  Then, obviously, $\delta (\eta_4^=)=(n/2)$, and $\sigma_2(\eta_4^=) =n $. See Fig. 8. 
\end{proof}
\section{Concluding Remark}
\indent \indent Following the previous study in [6], we are interested in the fault-tolerance of graphs satifying the degree conditions given by Ore [3]. Since the 1-fault tolerance has been thoroughly studied in [6], we further explore the 2-fault tolerance for any graph $G$ with $\sigma_2(G) \ge n$, 
where $|G|=n$. This paper has given a conclusion that any $G$ with  $\sigma_2(G)\ge n$ and $\kappa (G)\ge 4$ must be 2-vertex-fault tolerant unless $G$ belongs to one of the nine graph families in Theorem 6. Whether a given graph $G$ under the same condition is 2-edge-fault tolerant, or 1-vertex and 1-edge-fault tolerant, and how the exceptional graph families would be, remain open  problems. 
 \\

\end{document}